\documentclass{ifacconf}
\usepackage{graphicx}      
\usepackage{natbib}        
\usepackage{amsmath}
\usepackage{amssymb}

\newtheorem{theorem}{Theorem}
\newtheorem{remark}{Remark}

\newtheorem{lemma}{Lemma}
\newtheorem{assumption}{Assumption}
\newtheorem{corollary}{Corollary}

\begin{document}
\begin{frontmatter}

\title{On the Fundamental Limit of the Stochastic Gradient Identification Algorithm Under Non-Persistent Excitation\thanksref{footnoteinfo}} 

\thanks[footnoteinfo]{\textcopyright~2026 the authors. This work has been accepted to IFAC World Congress 2026 for publication under a Creative Commons Licence CC-BY-NC-ND. This research was supported by the National Natural Science Foundation of China under Grant No. 12288201 and by the National Key R\&D Program of China under Grant No. 2024YFC3307201. $^\dagger$ These authors contributed equally to this work.}

\author[First,Second]{Senhan Yao$^\dagger$}
\author[First,Second]{Longxu Zhang$^\dagger$}

\address[First]{State Key Laboratory of Mathematical Sciences, AMSS, Chinese Academy of Sciences, Beijing 100190, China (e-mail: yaosenhan@amss.ac.cn, zhanglongxu@amss.ac.cn).}
\address[Second]{School of Mathematical Sciences, University of Chinese Academy of Sciences, Beijing 100049, China.}

\begin{abstract}                
Stochastic gradient (SG) methods are fundamental to system identification and machine learning, enabling online parameter estimation in large-scale and streaming-data settings. As a classical identification method, the SG algorithm has been extensively studied for decades. Under non-persistent excitation, the strongest currently available convergence result assumes that the condition number of the Fisher information matrix is \(O((\log r_n)^\alpha)\), where \(r_n = 1 + \sum_{i=1}^n \|\varphi_i\|^2\). Existing theory establishes strong consistency when \(\alpha \le 1/3\), whereas the same condition with \(\alpha > 1\) is insufficient to guarantee strong consistency. We prove that strong consistency holds throughout the range \(0 \le \alpha < 1\). The proof is based on a new algebraic framework that yields substantially sharper matrix norm bounds. This result nearly resolves the four-decade-old Chen--Guo conjecture by establishing strong consistency throughout the previously open range \(1/3 < \alpha < 1\).
\end{abstract}


\end{frontmatter}

\section{Introduction}

\textit{How can we develop mathematical models of real-world physical processes from noisy observations to support critical engineering tasks such as controller design, prediction, and fault diagnosis?} This question lies at the heart of \textit{system identification}; see \cite{LJUNG20101}. In the era of artificial intelligence, characterized by large-scale datasets and an increasing demand for online learning, \textit{stochastic gradient (SG)} algorithms have regained prominence. These methods enable iterative online processing of noisy data, provide convergence guarantees for parameter estimates, and underlie optimization algorithms such as Adam; see \cite{ADAM}.

SG methods are closely related to \textit{stochastic approximation (SA)}, a framework originating from the pioneering work of \cite{robbins1951stochastic}. Subsequent milestones include the stochastic optimization approach developed by \cite{kiefer1952stochastic}, the general convergence theorem of \cite{dvoretzky1956stochastic}, asymptotic analyses by \cite{chung1954stochastic} and \cite{sacks1958asymptotic}, and the ODE method of \cite{gladyshev1965stochastic}. Together, these developments laid a rigorous mathematical foundation for this class of stochastic iterative algorithms.

Applications of these ideas in system identification emerged in the 1960s. \cite{astrom1965numerical} established a systematic \textit{offline} identification framework based on maximum likelihood estimation and introduced the pivotal concept of \textit{persistent excitation (PE)}. Concurrently, \cite{sakrison1964application} pioneered the application of stochastic approximation to \textit{online} identification, shifting the focus from offline modeling to real-time estimation. In the 1970s, attention shifted toward closed-loop systems. The groundbreaking work of \cite{astrom1972self} initiated the theoretical analysis of stochastic adaptive control, while \cite{ljung1977analysis} extended the ODE method of \cite{gladyshev1965stochastic} into a unified tool for analyzing the asymptotic behavior of general stochastic algorithms.

Although SG algorithms are known to produce strongly consistent estimates under PE conditions (see \cite{ljung1977analysis, anderson1979strong, chen1981strong}), their convergence behavior in the absence of PE warrants careful investigation. One motivation for studying this issue comes from adaptive control (see \cite{goodwin1980discrete}), where \textit{parameter estimation and optimal control can become decoupled}. Rigorously enforcing PE ensures strong consistency but may not yield an optimal controller; see \cite{chen1984strong}. Conversely, pursuing an optimal controller can compromise parameter consistency; see \cite{becker1985adaptive}. A second practical challenge stems from the fact that \textit{many systems operate under insufficiently rich inputs}, a phenomenon reported or studied in systems such as quadrotor UAVs (see \cite{chowdhary2012experimental}), neural networks (see \cite{nar2019persistencyexcitationrobustnessneural}), and spherical parallel robots (see \cite{RAD2020104026}). In such scenarios, inadequate excitation, if not properly accounted for, can lead to severely biased parameter estimates, unacceptably slow convergence, or even outright divergence of the identification algorithm.

The studies of \cite{chen1982strong} and \cite{lai1982} made independent and seminal contributions to the systematic relaxation of the PE condition for strong consistency of least-squares (LS) estimators. \cite{chen1982strong} established that strong consistency holds if the condition number of the Fisher information matrix satisfies
\begin{equation}\label{eq:chen-condition-number-growth}
    \kappa \left(\sum_{i=1}^n \varphi_i \varphi_i^\top \right) = O(r_n^\delta) \quad \mathrm{a.s.}, \quad 0 \leq \delta < \frac{1}{2},
\end{equation}
where \(r_n = 1 + \sum_{i=1}^n \|\varphi_i\|^2\). Simultaneously, Lai and Wei, working within a general stochastic regression framework, showed that strong consistency holds if
\begin{equation}\label{eq:lai-wei-condition-number-growth}
    \kappa \left(\sum_{i=1}^n \varphi_i \varphi_i^\top \right) = o\left(\frac{r_n}{\log r_n} \right) \quad \mathrm{a.s.}
\end{equation}

Building on this line of work, the Chen--Guo approach provides a deeper understanding of the SG algorithm under relaxed excitation conditions. \cite{chen1985strong} characterized the convergence rate of the SG algorithm under non-persistent excitation. This approach relates strong consistency to the convergence of an instrumental matrix sequence through a \textit{sample-pathwise} analysis; see \cite{chen1985consistency}. Under mild noise conditions and without requiring the noise to be i.i.d.\ or to have bounded variance, \cite{chen1985strong, guo1993, chen1986limit} showed that if
\begin{equation}\label{eq:condition-number-growth}
    \kappa \left(\sum_{i=1}^n \varphi_i \varphi_i^\top \right) = O((\log r_n)^\alpha) \quad \mathrm{a.s.}, 
\end{equation}
then the following statements hold:
\begin{itemize}
    \item \(\alpha \leq 1/3\) ensures strong consistency.
    \item \(\alpha > 1\) does not guarantee strong consistency.
\end{itemize}
This stark dichotomy inevitably raises the following fundamental questions: \textit{What are the fundamental limits of SG identification algorithms? To what extent can we rely on the outputs of SG algorithms?} More broadly, this dichotomy challenges us to understand how much insight into a system can be gained from finite, noisy, and poorly structured observations in a world replete with uncertainty.

\cite{chen1986limit} conjectured that the previously established sufficient condition is inherently conservative. They further conjectured that strong consistency should hold for the entire range \(0 \leq \alpha \leq 1\), with \(\alpha = 1\) representing the ``critical excitation'' threshold for the algorithm. For decades, the validity of this conjecture---and thus the characterization of the algorithm's fundamental limits within the gap \(1/3 < \alpha \leq 1\)---has remained open, posing a significant challenge in the theoretical foundations of SG identification.

The main contribution of this paper is \textit{a near-complete resolution of this long-standing conjecture.} Our results bridge the theoretical gap as follows:
\begin{itemize}
    \item We prove that, for the entire range \(0 \leq \alpha < 1\), condition~\eqref{eq:condition-number-growth} \textit{indeed suffices} for the almost-sure convergence of the SG estimates to the true parameter value, thereby providing a near-complete characterization of the convergence regime of the SG algorithm.
    \item To establish this result, we develop a novel \textit{algebraic approach} that provides a more transparent and versatile framework for obtaining the requisite matrix norm bounds, overcoming previous analytical obstacles.
\end{itemize}

The remainder of this paper is organized as follows. Section~\ref{sec:chen-guo-approach} revisits the problem formulation and the Chen--Guo approach. Section~\ref{sec:yao-framework} develops our new algebraic approach for bounding matrix products. Section~\ref{sec:zhang-analysis} presents detailed integral estimates. Section~\ref{sec:main-theorem} proves the main theorem, and Section~\ref{sec:concluding-remarks} concludes the paper.

\textbf{Notation.} The logarithm \(\log\) denotes the natural logarithm with base \(e\). We write \(A_n \lesssim B_n\) to mean that there exists a positive constant \(C>0\), independent of \(n\), such that \(A_n \leq C B_n\). The Landau symbols \(O(\cdot)\) and \(\Theta(\cdot)\) have their standard asymptotic meanings. \(\|\cdot\|\) denotes the Euclidean norm for vectors and the spectral norm for matrices unless specified otherwise. \(\|\cdot\|_F\) denotes the Frobenius norm. For a real symmetric matrix \(A\), \(\lambda_{\min}(A)\) and \(\lambda_{\max}(A)\) denote its smallest and largest eigenvalues, respectively. For a nonsingular matrix \(A\), the condition number \(\kappa(A)\) is defined by \(\kappa(A) = \|A\| \|A^{-1}\|\).

\section{Problem Formulation}\label{sec:chen-guo-approach}

In this section, we formulate the problem and briefly review the Chen--Guo approach. Consider the following multiple-input multiple-output (MIMO) system:
\begin{equation}
\begin{aligned}
 y_n + P_1 y_{n-1} &+ \cdots + P_p y_{n-p} \\
&=  Q_1 u_{n-1} + \cdots + Q_q u_{n-q} + \varepsilon_n. \\
\end{aligned}
\end{equation}
Here \(y_n\) and \(u_n\) denote the \(d\)-dimensional output and the \(m\)-dimensional input, respectively, and \(P_i\), \(i=1,\ldots,p\), and \(Q_j\), \(j=1,\ldots,q\), are the unknown system matrices. The noise \(\varepsilon_n\) is \(d\)-dimensional and is generated by a martingale difference sequence \(\{w_n\}\) as
\begin{equation}
\varepsilon_n = w_n + R_1 w_{n-1} + \cdots + R_r w_{n-r},
\end{equation}
where \(R_k\), \(k=1,\ldots,r\), are unknown matrices. The sequence \(\{w_n\}\) satisfies
\begin{equation}
\mathbb{E}[w_n|\mathcal F_{n-1}]=0,\quad \forall n\ge 1,
\end{equation}
where \(\{\mathcal F_n\}\) is an increasing family of \(\sigma\)-algebras on the probability space \((\Omega,\mathcal F,\mathbb{P})\).

Let \(z\) be the backward shift operator, and define
\begin{equation}
    R(z) = I + R_1 z + \cdots + R_r z^r.
\end{equation}
Set
\begin{equation}
\theta^\top =
\begin{bmatrix}
    -P_1 & \cdots & -P_p &
    Q_1 & \cdots & Q_q &
    R_1 & \cdots & R_r
\end{bmatrix}.
\end{equation}

To distinguish the exact regression representation from the recursion
used by the algorithm, we introduce two regressors. First, if
the true innovations were available, the ideal regressor would be
\begin{equation}
    \begin{aligned}
        \psi_n^\top =
        \big[
        &y_n^\top, \cdots, y_{n-p+1}^\top,
        u_n^\top, \cdots, u_{n-q+1}^\top, \\
        &w_n^\top, \cdots, w_{n-r+1}^\top
        \big].
    \end{aligned}
\end{equation}
With this ideal regressor, the ARMAX model admits an exact linear
regression form
\begin{equation}
    y_{n+1} = \theta^\top \psi_n + w_{n+1}.
\end{equation}
Note that a linear regression model can be viewed as an ARMAX model with \(p=r=0\) and \(q=1\). However, the innovations \(w_n, w_{n-1}, \ldots, w_{n-r+1}\) are not
observable. In practice, they are replaced by the estimated residuals
\begin{equation}
    \hat{w}_{n} := y_{n} - \theta_{n-1}^\top \varphi_{n-1}.
\end{equation}
This yields the computable regressor
\begin{equation}
    \begin{aligned}
        \varphi_n^\top =
        \big[
        &y_n^\top, \cdots, y_{n-p+1}^\top,
        u_n^\top, \cdots, u_{n-q+1}^\top, \\
        &\hat{w}_n^\top, \cdots, \hat{w}_{n-r+1}^\top
        \big].
    \end{aligned}
\end{equation}
Thus, \(\varphi_n\) is defined recursively and is
\(\mathcal{F}_n\)-measurable. Note that, when \(r>0\),
\(\varphi_n\) is generally different from the ideal regressor
\(\psi_n\). Consequently, the identity
\(y_{n+1}=\theta^\top\varphi_n+w_{n+1}\) is not exact in general; rather,
\begin{equation}
    y_{n+1}
    =
    \theta^\top\varphi_n
    + w_{n+1}
    + \theta^\top(\psi_n-\varphi_n).
\end{equation}
The last term represents the perturbation caused by replacing the true
innovations with their estimated residuals.

The estimation problem is then addressed through the following
regression-type stochastic gradient recursion. Denote by \(\theta_n\)
the estimate of \(\theta\) at time \(n\). Given deterministic initial
values \(\theta_0\) and \(\varphi_0\), define
\begin{equation}
    \theta_{n+1}
    =
    \theta_n
    +
    \frac{\varphi_n}{r_n}
    \left(
        y_{n+1}^\top - \varphi_n^\top \theta_n
    \right),
\end{equation}
where
\begin{equation}\label{eq:definition-r-n}
    r_n = 1 + \sum_{i=1}^n \|\varphi_i\|^2,
    \qquad r_0 = 1.
\end{equation}

Chen and Guo introduced the following instrumental transition matrix:
\begin{equation}
    \Phi(n+1,i) = (I - A_n)\Phi(n,i), \quad n \geq i,
\end{equation}
\begin{equation}
    \Phi(i,i) = I,
\end{equation}
where
\begin{equation}
    A_n = \frac{\varphi_n \varphi_n^\top}{r_n}.
\end{equation}

Chen and Guo used the \textit{strictly positive real (SPR)} condition to bound the error introduced by the estimated residuals, thereby ensuring that the parameter error vanishes as the instrumental transition matrix converges to zero.

\begin{theorem}[Chen and Guo, 1985a]
    If \(r=0\), or if \(r>0\) and \(R(z) - \frac{1}{2}I\) is SPR, then \(\Phi(n,0) \to 0\) implies \(\theta_n \to \theta\).
\end{theorem}

Notably, under mild noise assumptions, this condition \(\Phi(n,0) \to 0\) is also necessary when \(r=0\).

\begin{assumption}[``Condition A'']
    \ 

    \begin{itemize}
        \item As $n \to \infty$, $\sum_{i=0}^{n} \frac{\varphi_i}{r_i} \varepsilon_{i+1}^\top$ converges to a finite limit $S$.
        \item There exist $c > 0$ and $\delta > 0$, possibly depending on \(\omega\), such that \(
\left\| S - \sum_{i=0}^{n-1} \frac{\varphi_i}{r_i} \varepsilon_{i+1}^\top \right\| \leq c r_n^{-\delta} \) for all \(n\).
    \end{itemize}
\end{assumption}
\begin{remark}
    ``Condition A'' is a sample-path condition controlling the accumulated noise effect and is weaker than many standard independence assumptions.
\end{remark}

\begin{theorem}[Chen and Guo, 1985b]
    Assume that \(r=0\) and the noise sequence \(\{\varepsilon_n\}\) satisfies ``Condition A'' along a sample path \(\omega \in \Omega\). Then, along this sample path, for any initial value \(\theta_0\), we have \(\theta_n \to \theta\) if and only if \(\Phi(n,0) \to 0\), and in this case, the convergence rate is
    \begin{equation}
        \|\theta_n - \theta\| = O(\|\Phi(n,0)\|^{\delta/(1+\delta)}),
    \end{equation}
    where \(\delta > 0\) may depend on the sample path \(\omega\).
\end{theorem}

Thus, establishing strong consistency of the SG algorithm reduces to proving that $\Phi(n,0) \to 0$, which amounts to bounding the norm of a deterministic matrix product.

\begin{assumption}\label{ass:w-n}
    The sequence \(\{w_n\}\) is adapted to \(\{\mathcal F_n\}\) and satisfies

    \begin{itemize}
        \item \(\mathbb{E}[w_n\mid \mathcal{F}_{n-1}] = 0\).
        \item \(\mathbb{E}[\|w_n\|^2\mid \mathcal{F}_{n-1}] \leq c_0 r_{n-1}^\eta\), where \(c_0 > 0\), \(0 \leq \eta \leq 1\).
    \end{itemize}
\end{assumption}
\begin{remark}
    The bound $r_{n-1}^\eta$ allows for potential growth in the conditional noise variance, with $\eta = 0$ corresponding to the bounded-variance case and $\eta > 0$ permitting variance that grows with the accumulated regressor energy $r_n$. Assumption~\ref{ass:w-n} can often be used to verify that ``Condition A'' holds.
\end{remark}

\begin{assumption}\label{ass:r-n}
    \(r_n \to \infty\), and \(r_n = O(r_{n-1})\).
\end{assumption}
\begin{remark}
    The condition \(r_n \to \infty\) ensures that sufficient information is available for identification. In real-world applications, many physical systems exhibit bounded input-output behavior due to physical constraints, actuator limits, and sensor ranges. The condition $r_n = O(r_{n-1})$ naturally arises in such scenarios, as it implies that the energy injected into the system cannot grow arbitrarily fast between consecutive time steps.
\end{remark}

Using analytical techniques, Chen and Guo significantly relaxed the excitation requirements for strong consistency of parameter estimates and provided quantitative convergence rates.

\begin{theorem}[Guo, 1993]\label{thm:Guo1993}
    If ``Condition A'', Assumption~\ref{ass:r-n}, and the condition-number growth condition
    \begin{equation}
        \kappa\left( \sum_{i=1}^n \varphi_i \varphi_i^\top \right) = O((\log r_n)^{1/3})
    \end{equation}
     hold along a sample path \(\omega \in \Omega\), then \(\Phi(n,0) \to 0\) along that path. Moreover, the convergence rate is
    \begin{equation}
        \|\theta_n - \theta\| = O((\log r_n)^{-\delta}),
    \end{equation}
    where \(\delta>0\) may depend on the sample path \(\omega\).
\end{theorem}
\begin{remark}
The exponent \(1/3\) in Theorem~\ref{thm:Guo1993} was the largest known sufficient threshold before the present work.
\end{remark}

Using a counterexample, Chen and Guo showed that, if the condition number grows at a super-logarithmic rate, one can construct examples for which the SG algorithm fails.

\begin{theorem}[Chen and Guo, 1986]
    Suppose that Assumption~\ref{ass:w-n} is satisfied almost surely. Then, for any \(\delta > 0\), there exists a sequence of random vectors \(\{\varphi_n\}\) satisfying Assumption~\ref{ass:r-n} and the condition-number growth condition
    \begin{equation}
        \kappa\left( \sum_{i=1}^n \varphi_i \varphi_i^\top \right) = O((\log r_n)^{1+\delta}), \quad \mathrm{a.s.}
    \end{equation}
    but \(\Phi(n,0) \nrightarrow 0\) a.s.
\end{theorem}

Chen and Guo posited that the established threshold for the condition-number growth rate might be conservative, suggesting that the critical exponent in the bound could be sharpened to the limiting logarithmic order \(O(\log r_n)\).

\section{Mathematical Framework}\label{sec:yao-framework}

In this section, we consider a sequence \(\{A_n\}\) of symmetric matrices such that each \(A_n\) satisfies \(0 \leq A_n \leq I\), has rank at most one, and admits a decomposition of the form
\begin{equation}
    A_n = \phi_n \phi_n^\top.
\end{equation}
Our goal is to estimate products involving the instrumental transition matrix (Theorem~\ref{thm:instrumental-matrix-product-bound}). We introduce an auxiliary sequence:
\begin{equation}
    x_{i+1} = (I- A_i)x_i, \quad i \geq k.
\end{equation}
From this, we obtain
\begin{equation}\label{eq:x-n-splitting-term-2}
        x_i - x_k = -\sum_{j=k}^{i-1} A_j x_j,
    \end{equation}
and since \(A_{i-1}^2 \leq A_{i-1}\), it follows that
\begin{equation}\label{eq:x-n-splitting-term}
\|x_i\|^2 \leq \|x_{i-1}\|^2 - \langle A_{i-1} x_{i-1}, x_{i-1}\rangle.
\end{equation}
Summing these inequalities over \(j\) yields
\begin{equation}
        \sum_{j=k}^{i-1} \|\phi_j^\top x_{j}\|^2 \leq \|x_k\|^2 - \|x_i\|^2.
\end{equation}

We now introduce a nonnegative real sequence \(\mu_n \geq 0\) referred to as the \textit{weights}, and define a weighted sum \(S_{ik}\) of the matrices \(A_j\) over the interval \([k,i)\):
\begin{equation}
    S_{ik} = \sum_{j = k}^{i-1} \mu_j A_j.
\end{equation}

\begin{remark}
    The design of \(S_{ik}\) represents one of the key ingredients in this framework. The weights \(\mu_j\) serve multiple purposes: they can compensate for non-uniform regressor magnitudes, emphasize periods of high information content, or discount older measurements in time-varying systems.
\end{remark}

We compute the quadratic form:
\begin{equation}
    \begin{aligned}
        x_k^\top S_{ik} x_k & = x_k^\top \left(\sum_{j=k}^{i-1} \mu_j A_j \right) x_k =  \sum_{j=k}^{i-1} \mu_j  x_k^\top A_j x_k\\
        & =\sum_{j=k}^{i-1} \mu_j  x_k^\top \phi_j \phi_j^\top x_k = \sum_{j=k}^{i-1} \mu_j \|\phi_j^\top x_k\|^2.
    \end{aligned}
\end{equation}

To analyze this quantity, define the vectors
\begin{equation}
    \alpha = \begin{bmatrix}
        \ldots& \phi_j^\top x_k& \ldots
    \end{bmatrix}^\top,
\end{equation}
\begin{equation}
    \beta = \begin{bmatrix}
        \ldots& \phi_j^\top x_j& \ldots
    \end{bmatrix}^\top.
\end{equation}
Note that
\begin{equation}
    \beta-\alpha = \begin{bmatrix}
        \ldots& \phi_j^\top (x_j - x_k)& \ldots
    \end{bmatrix}^\top.
\end{equation}
Multiplying both sides of (\ref{eq:x-n-splitting-term-2}) by \(\phi_j^\top\) gives
\begin{equation}
    (\beta-\alpha)_j  = \phi_j^\top x_j - \phi_j^\top x_k = -\sum_{l=k}^{j-1} (\phi_j^\top \phi_l) (\phi_l^\top x_l).
\end{equation}
Since \(\phi_l^\top x_l = \beta_l\), we define a strictly lower-triangular matrix \(C\) by
\begin{equation}
    C_{jl} = \phi_j^\top \phi_l, \quad k \leq l < j.
\end{equation}

By construction, we have the matrix identity
\begin{equation}
    \alpha = (I+C)\beta.
\end{equation}

\begin{remark}
   The matrix \(C\) encodes the \textit{intertemporal correlation structure} of the regressor sequence, i.e., it quantifies how much information each new regressor \(\phi_j\) shares with previous regressors \(\phi_l\), \(l < j\). The strictly lower-triangular structure reflects the causal nature of time: future regressors cannot affect past ones. The matrix identity reveals that the initial projection error \(\alpha_j\) equals the current projection error \(\beta_j\) plus a correction term that accounts for how much the state has evolved due to previous updates.
\end{remark}

Set \(\Lambda = \mathrm{diag}(\ldots,\sqrt{\mu_j},\ldots)\). Then
\begin{equation}\label{eq:estimate-quadratic-form}
    \begin{aligned}
        \lambda_{\min}(S_{ik}) \|x_k\|^2 &\leq x_k^\top S_{ik} x_k
         = \|\Lambda \alpha \|^2 \\ &= \|\Lambda (I+C) \beta\|^2  \leq \|\Lambda (I+C) \|^2 \|\beta\|^2 \\
        & \leq \|\Lambda (I+C) \|^2 (\|x_k\|^2 - \|x_i\|^2).
    \end{aligned}
\end{equation}

It is a standard fact that the operator norm of a matrix can be bounded by its Frobenius norm:
\begin{equation}
    \|A\| \leq \|A\|_F.
\end{equation}
Hence, we obtain the estimate:
\begin{equation}\label{eq:estimate-lambda-i-c}
    \begin{aligned}
        \|\Lambda (I+C) \| &  \leq \|\Lambda\| + \|\Lambda C\| \leq \|\Lambda\| + \|\Lambda C\|_F\\
        & = \sqrt{\max_{k \leq j < i} \mu_j} + \sqrt{\sum_{j = k}^{i-1} \mu_j \sum_{l=k}^{j-1} (\phi_j^\top \phi_l)^2}.
    \end{aligned}
\end{equation}
\begin{remark}
    Estimating \(\Lambda\) and \(\Lambda C\) separately is important because it distinguishes between two fundamentally different sources of ``complexity'' in the system: the \textit{magnitudes} of the weights (controlled by \(\|\Lambda\|\)) and the \textit{temporal correlation structure} (captured by \(\|\Lambda C\|\)). In applications, this separation allows one to control system behavior independently through the choice of the weights \(\mu_j\) and through regressor design, which affects the correlation structure \(C\).
\end{remark}
For convenience, define
\begin{equation}\label{eq:definition-b-jk}
    B_{jk} = \sum_{l=k}^{j-1} (\phi_j^\top \phi_l)^2.
\end{equation}
Combining (\ref{eq:estimate-quadratic-form}), (\ref{eq:estimate-lambda-i-c}), and (\ref{eq:definition-b-jk}), we obtain the following inequality relating \(\|x_k\|\) and \(\|x_i\|\):
\begin{equation}
    \begin{aligned}
        \lambda_{\min}&(S_{ik}) \|x_k\|^2  \\
        & \leq \left[\sqrt{\max_{k \leq j < i} \mu_j} + \sqrt{\sum_{j = k}^{i-1} \mu_j B_{jk}}  \right]^2 (\|x_k\|^2 - \|x_i\|^2).
    \end{aligned}
\end{equation}

Since \(x_k\) is arbitrary, we may bound the norm of the instrumental transition matrix as follows.

\begin{theorem}\label{thm:instrumental-matrix-product-bound}
With the above notation, the following inequality holds:
\begin{equation}
    \|\Phi(N,k)\|^2 \leq  1 - \frac{\lambda_{\min}(S_{Nk})}{\left(\sqrt{\max_{k \leq j < N} \mu_j} + \sqrt{\sum_{j = k}^{N-1} \mu_j B_{jk}} \right)^2 }.
\end{equation}
\end{theorem}

We use the standard fact that, for any sequence \(\{a_i\}\) with \(0\le a_i<1\),
\begin{equation}
\prod_{i=1}^{\infty} (1-a_i) = 0 \Leftrightarrow \sum_{i=1}^{\infty} a_i = \infty.
\end{equation}

In our new framework, this basic fact yields an important corollary.

\begin{corollary}\label{cor:lambda-mu-b}
With the above notation, we have \(\Phi(n,0) \to 0\) if
\begin{equation}\label{eq:convergence-criterion-yao}
    \sum_{k =1}^\infty \frac{\lambda_{\min}(S_{t_k t_{k-1}})}{\left(\sqrt{\max_{t_{k-1} \leq j < t_k} \mu_j} + \sqrt{\sum_{j = t_{k-1}}^{t_k-1} \mu_j B_{jt_{k-1}}} \right)^2 }  = \infty,
\end{equation}
where \(\{t_k\}\) is a strictly increasing sequence of natural numbers tending to infinity, i.e., \(t_k \to \infty\).
\end{corollary}
\begin{remark}
This result relies on a ``time-scale rescaling'' technique. The idea is to partition the time axis into intervals $[t_{k-1}, t_k)$ and analyze the system over these aggregated blocks. This approach is necessary because, if the matrices were considered one at a time, i.e., over intervals of length one, the minimum eigenvalue $\lambda_{\min}(S_{t_k t_{k-1}})$ would be zero due to rank deficiency, causing the criterion to fail.
\end{remark}

The weight sequence \(\mu_j\) can be chosen as a \textit{streaming statistic} of the regressors \(\{\phi_j\}\), in the sense that each \(\mu_j\) is computed from \(\phi_j\) and \(\mu_{j-1}\). This property is particularly useful for online computation.

\section{Convergence Analysis}\label{sec:zhang-analysis}

In this section, we set \(\phi_j=\varphi_j/\sqrt{r_j}\), which yields the following estimate by the Cauchy--Schwarz inequality:
\begin{equation}
    \begin{aligned}
        B_{jk} & = \sum_{l=k}^{j-1} (\phi_j^\top \phi_l)^2 = \sum_{l=k}^{j-1}  \frac{(\varphi_j^\top \varphi_l)^2}{r_j r_l}  \\
        & \leq  \sum_{l=k}^{j-1} \frac{\|\varphi_j\|^2}{r_j}\frac{\|\varphi_l\|^2}{r_l} = \frac{\|\varphi_j\|^2}{r_j}  \sum_{l=k}^{j-1}  \frac{r_l - r_{l-1}}{r_l} \\
        & \leq \frac{\|\varphi_j\|^2}{r_j}  \int_{r_{k-1}}^{r_{j-1}} \frac{dx}{x}  =  \frac{r_j - r_{j-1}}{r_j}  \log \frac{r_{j-1}}{r_{k-1}}.
    \end{aligned}
\end{equation}
After substituting this estimate into the denominator of \eqref{eq:convergence-criterion-yao} and choosing \(\mu_j = r_j\), we obtain
\begin{equation}
\begin{aligned}
&\sum_{j = t_{k-1}}^{t_k-1} \mu_j B_{jt_{k-1}} \leq \sum_{j = t_{k-1}}^{t_k-1} (r_j - r_{j-1}) \log \frac{r_{j-1}}{r_{t_{k-1}-1}} \\
&\leq \int_{r_{t_{k-1}-1}}^{r_{t_{k}-1}} \log x \, dx - (r_{t_k -1} - r_{t_{k-1}-1}) \log r_{t_{k-1}-1} \\
&= r_{t_{k}-1} \log r_{t_{k}-1} - r_{t_{k}-1} - r_{t_{k-1}-1} \log r_{t_{k-1}-1} \\
&\quad + r_{t_{k-1}-1}  - (r_{t_k -1} - r_{t_{k-1}-1}) \log r_{t_{k-1}-1} \\
&= r_{t_{k}-1} \left( \log r_{t_{k}-1} - \log r_{t_{k-1}-1} \right) - (r_{t_{k}-1} - r_{t_{k-1}-1}).\\
\end{aligned}
\end{equation}

This yields an explicit sufficient criterion for \(\Phi(n,0)\to 0\).
\begin{corollary}\label{cor:zhang-criterion}
With the above notation, we have \(\Phi(n,0) \to 0\) if
\begin{equation}\label{eq:zhang-criterion}
    \sum_{k =1}^\infty \frac{\lambda_{\min}(S_{t_k t_{k-1}})}{D_k} = \infty,
\end{equation}
where \(\{t_k\}\) is defined as in Corollary~\ref{cor:lambda-mu-b}, \(S_{t_k t_{k-1}} = \sum_{i=t_{k-1}}^{t_k - 1} \varphi_i \varphi_i^\top\), and \(D_k = r_{t_{k}-1} \left( \log r_{t_{k}-1} - \log r_{t_{k-1}-1} \right) + r_{t_{k-1}-1}\).
\end{corollary}

\section{Main Theorem}\label{sec:main-theorem}

The following lemma is analogous to the ``time-inverse function'' introduced by Chen and Guo, but is formulated in a way that simplifies the subsequent analysis.

\begin{lemma}\label{lem:choose-t-k}
Assume that Assumption~\ref{ass:r-n} holds. Then there exist a strictly increasing sequence of natural numbers \(\{t_k\}\) and a constant \(L>1\) such that
\begin{equation}
    \frac{k}{L} < \frac{r_{t_k}}{r_{t_{k-1}}} < Lk.
\end{equation}
\end{lemma}
\begin{pf}
Define \(t_k = \min \{j : r_j \geq k!\}\). By Assumption~\ref{ass:r-n}, there exists \(L>1\) such that \(r_n \leq Lr_{n-1}\) for all sufficiently large \(n\). Since \(t_k \to \infty\) and \(k+1>L\) eventually, after discarding finitely many initial terms and reindexing, we may take \(\{t_k\}\) to be strictly increasing. Moreover,
\begin{equation}
    k! \leq r_{t_k} \leq L r_{t_k - 1} < Lk!,
    \quad
    (k-1)! \leq r_{t_{k-1}} < L(k-1)!.
\end{equation}
Combining these yields the desired inequality.
\end{pf}

We now prove the main theorem.

\begin{theorem}
    On any sample path on which Assumption~\ref{ass:r-n} holds and \(S_n = \sum_{i=1}^n \varphi_i \varphi_i^\top\) is eventually nonsingular, if
    \begin{equation}\label{eq:condition-number-growth-assumption-thm6}
        \kappa\left( \sum_{i=1}^n \varphi_i \varphi_i^\top \right) = O((\log r_n)^{\alpha}),
    \end{equation}
    where \(0 \leq \alpha < 1\), then \(\Phi(n,0) \to 0\).
\end{theorem}
\begin{pf}
    For convenience, set
    \begin{equation}
        S_n = \sum_{i=1}^n \varphi_i \varphi_i^\top.
    \end{equation}
    The block information matrix in \eqref{eq:zhang-criterion}, \(S_{t_k t_{k-1}}\), can be expressed as the difference of two cumulative matrices:
    \begin{equation}
        S_{t_k t_{k-1}} = S_{t_k-1} - S_{t_{k-1}-1}.
    \end{equation}
    We use Weyl's inequality for the eigenvalues of a sum of Hermitian matrices, which states that
    \begin{equation}
        \lambda_{\min}(S_{t_k t_{k-1}}) \geq \lambda_{\min}(S_{t_k-1}) - \lambda_{\max}(S_{t_{k-1}-1}).
    \end{equation}
    We next derive a lower bound for the first term \(\lambda_{\min}(S_{t_k-1})\) and an upper bound for the second term \(\lambda_{\max}(S_{t_{k-1}-1})\). By \eqref{eq:condition-number-growth-assumption-thm6}, for the fixed sample path under consideration, there exists \(M>0\) such that
    \begin{equation}
        \begin{aligned}
            \lambda_{\min}(S_{t_k-1}) &\geq \frac{\lambda_{\max}(S_{t_k-1})}{M(\log r_{t_k-1})^{\alpha}} \\
            & \gtrsim \frac{\mathrm{tr}(S_{t_k-1})}{M(\log r_{t_k-1})^{\alpha}} = \frac{r_{t_k-1} - 1}{M(\log r_{t_k-1})^{\alpha}}.
        \end{aligned}
    \end{equation}
    The upper bound follows from the monotonicity of \(r_n\):
    \begin{equation}
        \lambda_{\max}(S_{t_{k-1}-1}) \leq \mathrm{tr}(S_{t_{k-1}-1}) = r_{t_{k-1}-1} -1.
    \end{equation}
    Substituting these bounds into Weyl's inequality gives the following lower bound:
    \begin{equation}
        \lambda_{\min}(S_{t_k t_{k-1}}) \gtrsim \frac{r_{t_k-1} - 1}{M(\log r_{t_k-1})^{\alpha}} - ( r_{t_{k-1}-1} -1).
    \end{equation}
    Choose \(t_k\) as in Lemma~\ref{lem:choose-t-k}. By Stirling's formula, we have
    \begin{equation}
        \log r_{t_k-1} \leq \log (Lk!) \lesssim k\log k.
    \end{equation}
    This yields a lower bound for \(\lambda_{\min}(S_{t_k t_{k-1}})\) in terms of \(k\):
    \begin{equation}
        \begin{aligned}
            \lambda_{\min}(S_{t_k t_{k-1}}) & \gtrsim \frac{k!}{(k \log k)^\alpha} - (k-1)! \\
            & = (k-1)! \left( \frac{k^{1-\alpha}}{(\log k)^\alpha} - 1\right)  \gtrsim \frac{k!}{(k \log k)^\alpha}.
        \end{aligned}
    \end{equation}
    Moreover, note that
    \begin{equation}
        \begin{aligned}
            D_k & \leq  r_{t_{k}-1} \left( \log r_{t_{k}-1} - \log r_{t_{k-1}-1} + 1\right) \\
            & \lesssim r_{t_{k}} \log \frac{r_{t_{k}}}{r_{t_{k-1}}} < r_{t_k} \log (Lk)  \lesssim k! \log k.
        \end{aligned}
    \end{equation}
    To estimate the summand in (\ref{eq:zhang-criterion}), we obtain
    \begin{equation}\label{eq:comparison_lambda_D}
    \frac{\lambda_{\min}(S_{t_k t_{k-1}})}{D_k} \gtrsim \frac{k!}{(k \log k)^\alpha}\frac{1}{k! \log k} = \frac{1}{k^{\alpha} (\log k)^{1+\alpha}}.
    \end{equation}
    Since \(0\le \alpha<1\), the series whose general term is given by the right-hand side of \eqref{eq:comparison_lambda_D} diverges. Hence, by the comparison test, the series in \eqref{eq:zhang-criterion} diverges. It follows from Corollary~\ref{cor:zhang-criterion} that \(\Phi(n,0) \to 0\).
\end{pf}

\section{Conclusions}\label{sec:concluding-remarks}

This paper has nearly resolved a long-standing conjecture of Chen and Guo by demonstrating that the stochastic gradient algorithm achieves strong consistency even under non-persistent excitation, provided that the condition number of the Fisher information matrix is \(O((\log r_n)^\alpha)\) for some \(0 \leq \alpha < 1\). Our work broadens the known sufficient range from \(\alpha \leq 1/3\) to \(0 \leq \alpha < 1\) and identifies \(\alpha = 1\) as the remaining critical boundary. By introducing a novel algebraic framework, we have provided sharper matrix bounds and a more transparent proof. This advance deepens the theoretical understanding of stochastic gradient methods and highlights the remaining challenge of settling the boundary case \(\alpha = 1\). Future work may extend this framework to other step-size rules and investigate its implications for deep learning and adaptive control scenarios.

\begin{ack}
The authors thank Professor Lei Guo (AMSS) for raising the question addressed in this paper and for his valuable guidance. The authors also thank Yujing Liu (AMSS), Cheng Zhao (AMSS), Xin Zheng (AMSS), and the anonymous reviewers for helpful comments.
\end{ack}

\bibliography{ifacconf}

\appendix
\end{document}